\documentclass[10pt]{article}

\usepackage[final]{graphics}
\usepackage{graphicx}          	 
\usepackage{bm}                 
\usepackage{amsmath}             
\usepackage{amssymb}
\usepackage{amsfonts}             
\usepackage{verbatim}           
\usepackage{amsthm}             

\usepackage{mathtools}
\usepackage{relsize}

\numberwithin{equation}{section}	
\usepackage{hyperref}
\usepackage{xcolor}
\theoremstyle{plain}             

\newtheorem{theorem}{Theorem}[section]

\theoremstyle{definition}

\newtheorem{remark}[theorem]{Remark}

\input epsf
\def\protectbold#1{\protect{\boldmath{$#1$}}}

\def\sign{{\rm sign}}

\def\bigO{{\cal O}}

\def\RR{{\mathbb R}}
\def\CC{{\mathbb C}}
\def\ZZ{{\mathbb Z}}
\def\calC{{\cal C}}
\def\wt{\widetilde}

\def\phase{{\rm ph}}
\def\dsp#1{\displaystyle#1}
\def\protectbold#1{\protect{\boldmath{$#1$}}}

\begin{document}

\title{Asymptotic expansions of Kummer hypergeometric functions for large values of the parameters}
 \author{
Nico M. Temme\\
IAA, 1825 BD 25, Alkmaar, The Netherlands. \\
Former address: Centrum Wiskunde \& Informatica (CWI), \\
Science Park 123, 1098 XG Amsterdam, The Netherlands. \\
Email: Nico.Temme@cwi.nl
}

\date{\  }

\maketitle
\begin{abstract}
\noindent
New asymptotic expansions are derived of the Kummer functions $M(a,b,z)$ and $U(a,b+1,z)$ for large positive values of $a$ and $b$, with $z$ fixed. For both functions we consider $b/a\le 1$ and $b/a\ge 1$, with special attention for the case $a\sim b$. We use a uniform method to handle all cases of these parameters.
\end{abstract}

{\small
\noindent
{\bf Keywords} Asymptotic expansions; 
Kummer functions; Confluent hypergeometric functions\\
{\bf AMS Classification} Primary 41A60; Secondary 33C15}

\section{Introduction}\label{sec:intro}
Many asymptotic expansions of the Kummer functions (or confluent hypergeometric functions)  $M(a,b,z)$ and $U(a,b,z)$ are available in the literature. With the results of this paper we fill a gap regarding the case of large positive parameters $a$ and $b$, with real or complex argument  $z$  fixed or bounded. 

For $b\to\infty$, with $\vert z\vert \ll b$ and $a\ll b$, we can use the defining convergent power series given in \eqref{eq:appB01}, which has an asymptotic character. 
An asymptotic expansion in negative powers of $b$ can be found in \S13.8(i) of \cite{Olde:2010:CHF}, together with other asymptotic forms. We can also refer to  \cite[Chapter~10]{Temme:2015:AMI}, where several expansions of the Kummer functions for large $a$ or $b$ are considered.  Usually the available asymptotic relations  are in terms of the argument $z$ in combination with one or both parameters.  

In the present paper  we derive new asymptotic expansions of the Kummer functions $M(a,b,z)$ and $U(a,b+1,z)$ for large values of $a$ and $b$, with $z$ fixed. Special attention is required when $a\sim b$, in which case we derive expansions that are uniformly valid when the ratio $a/b$ approaches~1. We give new results for the following four cases, which are not considered earlier in the literature:
\begin{enumerate}
\item\quad  $M(a,b,z),\quad \quad \ b\ge a$; \quad \quad\S\ref{sec:Mbgea}; expansion \eqref{eq:Mbgea15}.
\item \quad $M(a,b,z),\quad \quad \ b\le a$; \quad \quad\S\ref{sec:Mblea}; expansion \eqref{eq:Mblea15}.
\item \quad $U(a,b+1,z),\quad b\ge a$; \quad \quad\S\ref{sec:Ubgea}; expansion \eqref{eq:Ubgea14}.
\item \quad $U(a,b+1,z),\quad b\le a$; \quad \quad \S\ref{sec:Ublea}. expansion \eqref{eq:Ublea12}.
\end{enumerate}
Throughout the paper we assume that both $a$ and $b$ are large, with $z=\bigO(1)$ and for the $U$-functions $\Re z > 0$. When $a$ or $b$ are of order $\bigO(1)$, the existing literature gives sufficient information.

For the asymptotics we use a rather simple uniform method to derive the large-$w$ asymptotic expansion of the Laplace-type integral
\begin{equation}\label{eq:intro01}
F_\lambda(w)=\frac1{\Gamma(\lambda)}\int_0^\infty
s^{\lambda-1}e^{-w\,s}f(s)\,ds,
\end{equation}
which expansion is uniformly valid  with respect to $\lambda\ge0$. A similar contour integral is also used.  We summarise this method in Appendix~A, using details of \cite[Chapter~25]{Temme:2015:AMI}. In Appendix~B we cite the most relevant formulas of the Kummer functions used in this paper.

\section{\protectbold{M(a,b,z),\ b \ge a}}\label{sec:Mbgea}
In this section we use the notation and condition
\begin{equation}\label{eq:Mbgea01}
\lambda = b-a,\quad \mu=\frac{\lambda}{a}=\frac{b-a}{a},\quad z\in\CC,\quad \vert z\vert \le z_0,
\end{equation}
where $z_0$ is a fixed positive number.
We use the Kummer relation for the $M$-function in \eqref{eq:appB07} together with \eqref{eq:appB02}. This gives
\begin{equation}\label{eq:Mbgea02}
M(a,b,z)=\frac{\Gamma(b)e^{z}}{\Gamma(a)\Gamma(\lambda)}\int_0^1 e^{-zt} e^{-a\phi(t)}\,\frac{dt}{t(1-t)},
\end{equation}
where
\begin{equation}\label{eq:Mbgea03}
\phi(t)=-\ln (1-t)-\mu\ln t.
\end{equation}
The saddle point $t_0$ follows from the zero of $\phi^\prime(t)$. We have
\begin{equation}\label{eq:Mbgea04}
\phi^\prime(t)=\frac{t(1+\mu)-\mu}{t(1-t)}\quad  \Longrightarrow \quad t_0=\frac{\mu}{1+\mu}.
\end{equation}
When the saddle point is properly inside the interval $[0,1]$ we can use the standard method for obtaining an asymptotic expansion by using the substitution $\phi(t)-\phi(t_0)=\frac12w^2$, $\sign(w)=\sign(t-t_0)$. However, when $t_0\to0$, that is, when $b\downarrow a$,
the standard method is no longer applicable, and we use a uniform method in which $b=a$ can be used.

The uniform method is based on a transformation of the integral in \eqref{eq:Mbgea02} into the standard form in \eqref{eq:intro01} by writing
\begin{equation}\label{eq:Mbgea05}
\phi(t)-\phi(t_0)=\psi(s)-\psi(s_0),\quad \sign(t-t_0)=\sign(s-s_0),
\end{equation}
where
\begin{equation}\label{eq:Mbgea06}
\psi(s)=s-\mu\ln s,\quad s_0=\mu;
\end{equation}
$s_0$ is the zero of $\psi^\prime(s)=(s-\mu)/s$.

\begin{figure}[tb]
\vspace*{0.8cm}
\begin{center}
\begin{minipage}{5cm}
   \includegraphics[width=5cm]{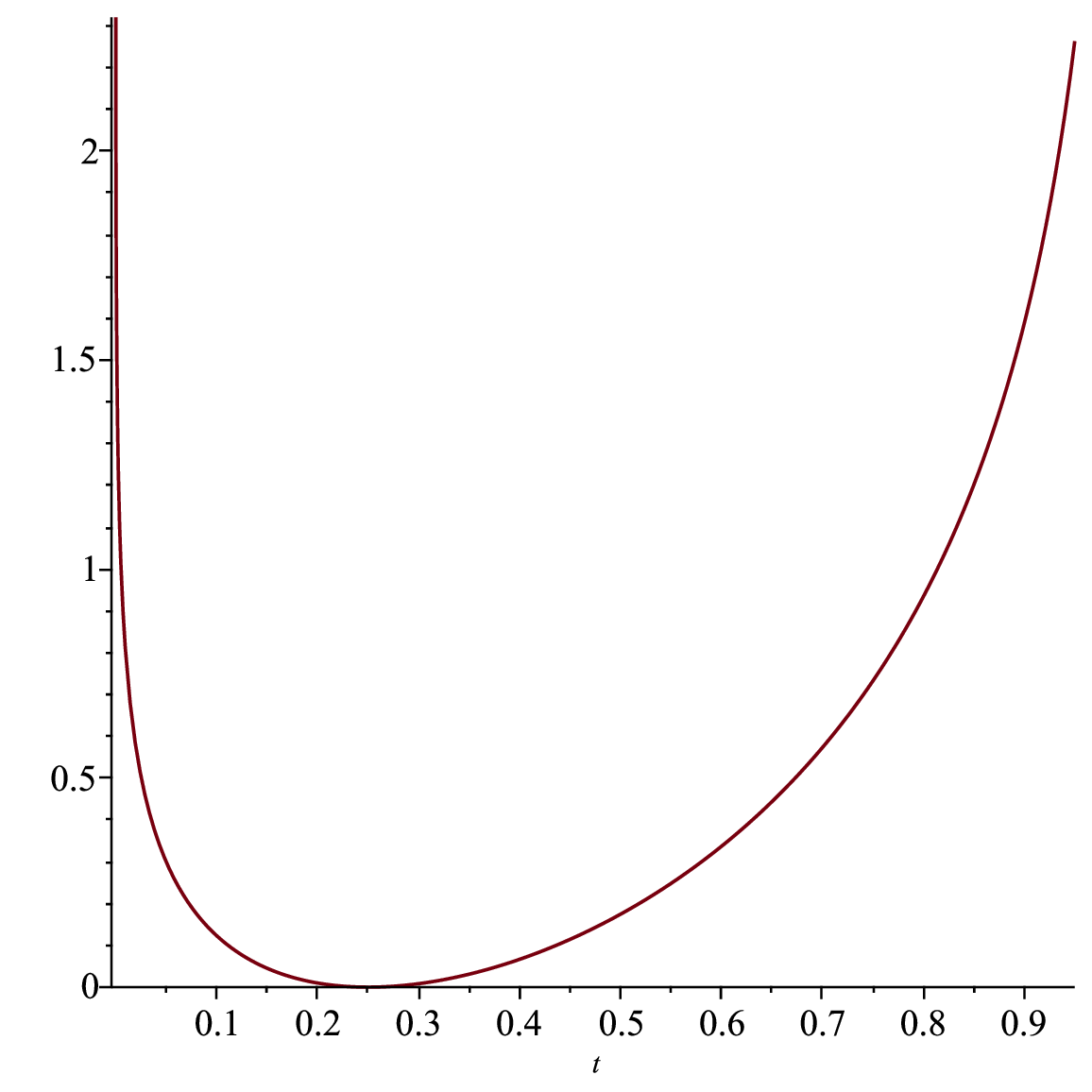} 
\end{minipage}
\hspace*{1cm}
\begin{minipage}{5cm}
   \includegraphics[width=5cm]{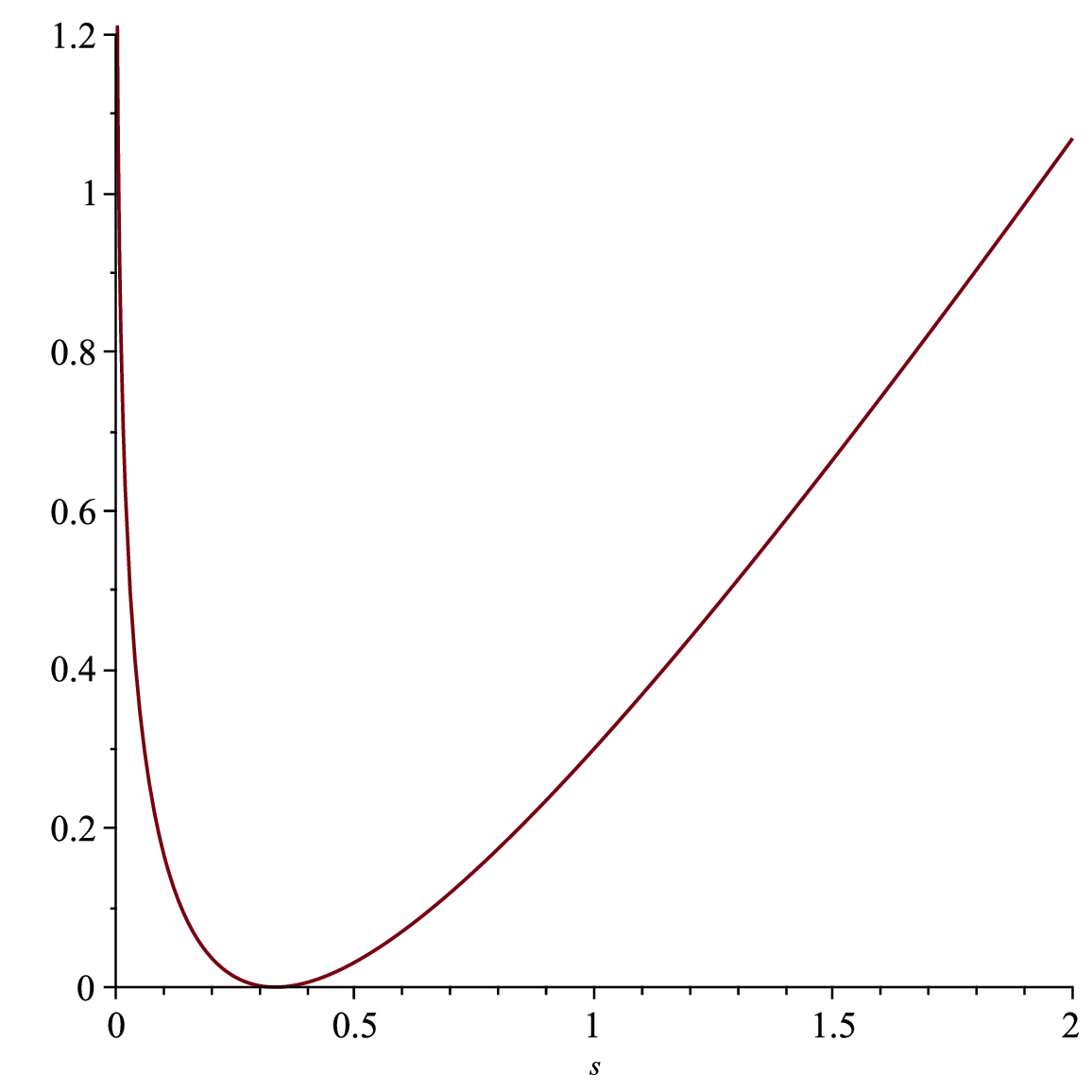} 
\end{minipage}
\end{center}
\caption{\small
Curves of the functions $\phi(t)-\phi(t_0)$ (left) and $\psi(s)-\psi(s_0)$ (right) that we use in the transformation in \eqref{eq:Mbgea05}, displayed for $\mu=\frac13$. }
\label{fig:fig01}
\end{figure}

In Figure~\ref{fig:fig01} we show the curves of the functions $\phi(t)-\phi(t_0)$ (left) and $\psi(s)-\psi(s_0)$ (right) that we use in the transformation in \eqref{eq:Mbgea05}; we use  $\mu=\frac13$. The convex curves touch the real axes at  $t_0=\frac{\mu}{1+\mu}=\frac14$ and $s_0=\mu=\frac13$. The condition $\sign(t-t_0)=\sign(s-s_0)$ means that the function values at the left of $t_0$ and $s_0$ correspond to each other, and the same holds true for those at the right of these points. Clearly, in this way, the transformation is one-to-one for $t\in(0,1)$ and $s >0$.

The transformation gives
\begin{equation}\label{eq:Mbgea07}
M(a,b,z)=\frac{\Gamma(b)}{\Gamma(a)}e^{z-aA}F_\lambda(a),\quad F_\lambda(a)=\frac{1}{\Gamma(\lambda)}\int_0^\infty e^{-as} s^ {\lambda-1} f(s)\,ds,
\end{equation}
where 
\begin{equation}\label{eq:Mbgea08}
A=\phi(t_0)-\psi(s_0)=(1+\mu)\ln(1+\mu)-\mu, \quad  f(s)=\frac{e^{-zt}}{1+\mu}\frac{s-\mu}{t-t_0},
\end{equation}
because
\begin{equation}\label{eq:Mbgea09}
f(s)=e^{-zt}\frac{s}{t(1-t)}\frac{dt}{ds},\quad \frac{dt}{ds}=\frac{\psi^\prime(s)}{\phi^\prime(t)}=\frac{s-\mu}{s}\frac{t(1-t)}{t(1+\mu)-\mu}.
\end{equation}

Using the expansion given in \eqref{eq:appA07}, we obtain
\begin{equation}\label{eq:Mbgea10}
M(a,b,z)\sim e^{z-aA}\frac{\Gamma(b)}{\Gamma(a)} a^{-\lambda}\sum_{n=0}^\infty
\frac{f_n(\mu)}{a^n},\quad a\to\infty,\quad b\ge a.
\end{equation}
The coefficients $f_n(\mu)$ are linear combinations of the derivatives of $f(s)$ at the saddle point $s=\mu$.

To find $f_0(\mu)$ we observe that in the definition of $f(s)$,
see \eqref{eq:Mbgea08} and \eqref{eq:Mbgea09}, we need the derivative $dt/ds$  at $s=\mu$. Because $s=\mu$ corresponds with $t=t_0$, we need to evaluate $dt/ds$  by using   l'H{\^o}pital's rule. We have
\begin{equation}\label{eq:Mbgea11}
\left.\frac{dt}{ds}\right\vert_{s=s_0}=
\frac{\psi^{\prime\prime}(s_0)}{\phi^{\prime\prime}(t_0)\left.\frac{dt}{ds}\right\vert_{s=s_0}}.
\end{equation}
This gives
\begin{equation}\label{eq:Mbgea12}
\left(\left.\frac{dt}{ds}\right\vert_{s=s_0}\right)^2=
\frac{\psi^{\prime\prime}(s_0)}{\phi^{\prime\prime}(t_0)}=\frac{1}{(1+\mu)^3}\quad\Longrightarrow \quad f_0(\mu)=e^{-zt_0}\sqrt{1+\mu}.
\end{equation}
We take the coefficient $f_0(\mu)$ in front of the expansion and write 
\begin{equation}\label{eq:Mbgea13}
M(a,b,z)\sim e^{z-aA}\frac{\Gamma(b)}{\Gamma(a)}a^{-\lambda}f_0(\mu)\sum_{n=0}^\infty \frac{\wt{f}_n(\mu)}{a^n},  \quad
\wt{f}_n(\mu)=\frac{f_n(\mu)}{f_0(\mu)}, \quad a,b\to\infty,\quad b\ge a.
\end{equation}
We evaluate the front factors by using the definition of $A$ in \eqref{eq:Mbgea08} and the scaled gamma functions defined in \eqref{eq:appB10}, and obtain
\begin{equation}\label{eq:Mbgea14}
e^{z-aA}\frac{\Gamma(b)}{\Gamma(a)}a^{-\lambda}f_0(\mu)=e^{z/(1+\mu)}\frac{\Gamma^*(b)}{\Gamma^*(a)}.
\end{equation}
This gives the final result
\begin{equation}\label{eq:Mbgea15}
M(a,b,z)\sim e^{z/(1+\mu)}\frac{\Gamma^*(b)}{\Gamma^*(a)}\sum_{n=0}^\infty \frac{\wt{f}_n(\mu)}{a^n}, \quad a,b\to\infty,\quad b\ge a.
\end{equation}
If we wish we can expand the ratio of scaled gamma functions in front of this expansion in powers of~$a^{-1}$,  using $b=a(1+\mu)$ (see \cite[\S6.5]{Temme:2015:AMI}).

The first few coefficients of this expansion are ${\wt f}_0(\mu)=1$,
\begin{equation}\label{eq:Mbgea16}
\begin{array}{r@{\,}c@{\,}l}
{\wt f}_1(\mu)&=&\dsp{\frac{\mu\left((\mu+1)^2+6z^2\right)}{12(\mu+1)^3},}\\[8pt]
{\wt f}_2(\mu)&=&\dsp{\frac{\mu\left(\mu(\mu+1)^4+12(\mu-12)(\mu+1)^2z^2+96(\mu^2-1)z^3+36\mu z^4\right)}{288(\mu+1)^6}.}
\end{array}
\end{equation}
These follow from the scheme given in Appendix~A. For the analytical  evaluation of these coefficients we refer to \S\ref{sec:num}, where also numerical details of the performance of the expansion are given.

\begin{remark}\label{rem:rem01}
To obtain a qualitative bound of the remainder $E_K$ of the expansion shown in  \eqref{eq:appA06}, we observe that
the function  $f(s)$ defined in \eqref{eq:Mbgea08} behaves as $f(s)=\bigO(s)$ as $s\to\infty$, because $t\in[0,1]$ and $z$ is assumed to be fixed. From the representation  in terms of rational functions in \eqref{eq:appA13}, and because $R_n(\sigma,s,\mu)=\bigO(1/s)$ for large $s$\footnote{This follows from the first functions given in \eqref{eq:appA14} and induction with respect to $n$.} we conclude that $f_n(s)=\bigO(1)$ for large $s$. We infer that the remainder $E_K$ in the finite expansion in \eqref{eq:appA06} for the present case is  $\bigO(1)$ with respect to the large parameter $a$. The rational functions are also bounded functions as $\mu\to\infty$. 
\end{remark}

\subsection{Details about the transformation}\label{sec:Mbgea2}
We give details about the transformation used in  \eqref{eq:Mbgea05}, the singularities of the function $f(s)$, and the uniform character of the expansion for $\mu\ge0$.  

The nonlinear transformation \eqref{eq:Mbgea05} can be inverted by using the Lambert $W$ function that satisfies the equation 
\begin{equation}\label{eq:Mbgea17}
W(z)e^{W(z)}=z.
\end{equation}
See \cite{Corless:1996:LWF} for details. For a proper description of $W(z)$ for $z\in\RR$ and $z\in\CC$, several branches of this function have to be considered. Write $s=-\mu\sigma$. Then for $\mu>0$  the transformation \eqref{eq:Mbgea05} can be written in the form
\begin{equation}\label{eq:Mbgea18}
\sigma e^\sigma=-\frac{t}{\mu}(1-t)^{1/\mu}e^{A(\mu)/\mu},
\end{equation}
where $A(\mu)$ is given in  \eqref{eq:Mbgea08}. We need to solve  this equation for $\sigma<0$, with the condition $\sign(\sigma+1)=\sign(t_0-t)$. For $\sigma=-1$ and $t=t_0$ both functions in  \eqref{eq:Mbgea18} have the value~$-1/e$.

In \cite{Temme:1985:LTI} we have shown that an expansion as the one obtained in \eqref{eq:Mbgea15}  is uniformly valid with respect to $\mu\ge0$ when $f(s)$ can be bounded by an algebraic function. Also,  the singularities of $f(s)$ should be bounded away from the positive axis, and  the distance of the singularities from the saddle point $s_0=\mu$ is larger than $d\sqrt{\mu}$, for some $d>0$. The singularities of the present function $f(s)$ satisfy these conditions. In some other sections we cannot give an algebraic bound. 

We can find the singularities by observing that these are generated by the multivalued logarithmic term $-\ln(1-t)$ of $\phi(t)$. The derivative $dt/ds$, which is part of $f(s)$, has singularities for $t$-values $t_0 e^{2\pi ik}$, $k\in\ZZ\setminus \{0\}$ outside the standard domain of the  logarithm;
$dt/ds$ is well defined for $k=0$. 
 
The singularities in the $s$-plane follow from the equation
\begin{equation}\label{eq:Mbgea19}
-\ln\left((1-t_0)e^{2\pi ik}\right)-\mu\ln(t_0)-\phi(t_0)=s-\mu\ln(s)-\mu+\mu\ln\mu,\quad k\ne0,
\end{equation}
or
\begin{equation}\label{eq:Mbgea20}
-2\pi ik=s-\mu\ln(s)-\mu+\mu\ln\mu,\quad k\ne0.
\end{equation}
There is no need to consider the logarithm $-\mu\ln t$ in the transformation, because we have 
chosen the logarithm in $\psi(s)$ with the same pre factor $\mu$. 
This gives an analytic relation between $t$ and $s$ at the origins.

The solutions $s_k(\mu)$ of \eqref{eq:Mbgea20} with $k=\pm1$ are closest to the domain of integration. We have
\begin{equation}\label{eq:Mbgea21}
s_{\pm1}(0)=\mp2\pi i,\quad s_{\pm1}(\mu)\sim \mu+2\sqrt{\pi\mu}\,e^{\mp\frac14\pi i},\quad \mu\to\infty.
\end{equation}
A graph given in  \cite{Temme:1985:LTI}  shows that indeed $\vert s_{\pm1}(\mu)-\mu\vert \ge d\sqrt{\mu}$ for some $d>0$.
For the loop integrals in the $s$-plane in later sections it is good to know that there are no singularities in the left half plane $\Re s <0$.

\section{\protectbold{M(a,b,z),\ b \le a}}\label{sec:Mblea}
In this section we use the notation and conditions

\begin{equation}\label{eq:Mblea01}
\lambda = a-b,\quad \mu=\frac{\lambda}{b}=\frac{a-b}{b},\quad \mu\le \mu_0,\quad z\in\CC,\quad \vert z\vert \le z_0,
\end{equation}
where $\mu_0$ and $z_0$ are  fixed positive numbers. The condition on $\mu$ means that $a/(1+\mu_0)\le b \le a$, and that, say, $b=o(a)$ is not allowed.

We use the integral representation given in \eqref{eq:appB03} and write it in the form
\begin{equation}\label{eq:Mblea02}
M(a,b,z)=\frac{\Gamma(b)\Gamma(1+\lambda)}{\Gamma(a)}\frac{1}{2\pi i} \int _0^{(1+)}e^{zt}e^{b\phi(t)}\frac{dt}{t(t-1)},
\end{equation}
where
\begin{equation}\label{eq:Mblea03}
\phi(t)=(1+\mu)\ln t-\mu\ln(t-1).\quad \
\end{equation}

\begin{figure}[tb]
\vspace*{0.8cm}
\begin{center}
\begin{minipage}{4cm}
   \includegraphics[width=5cm]{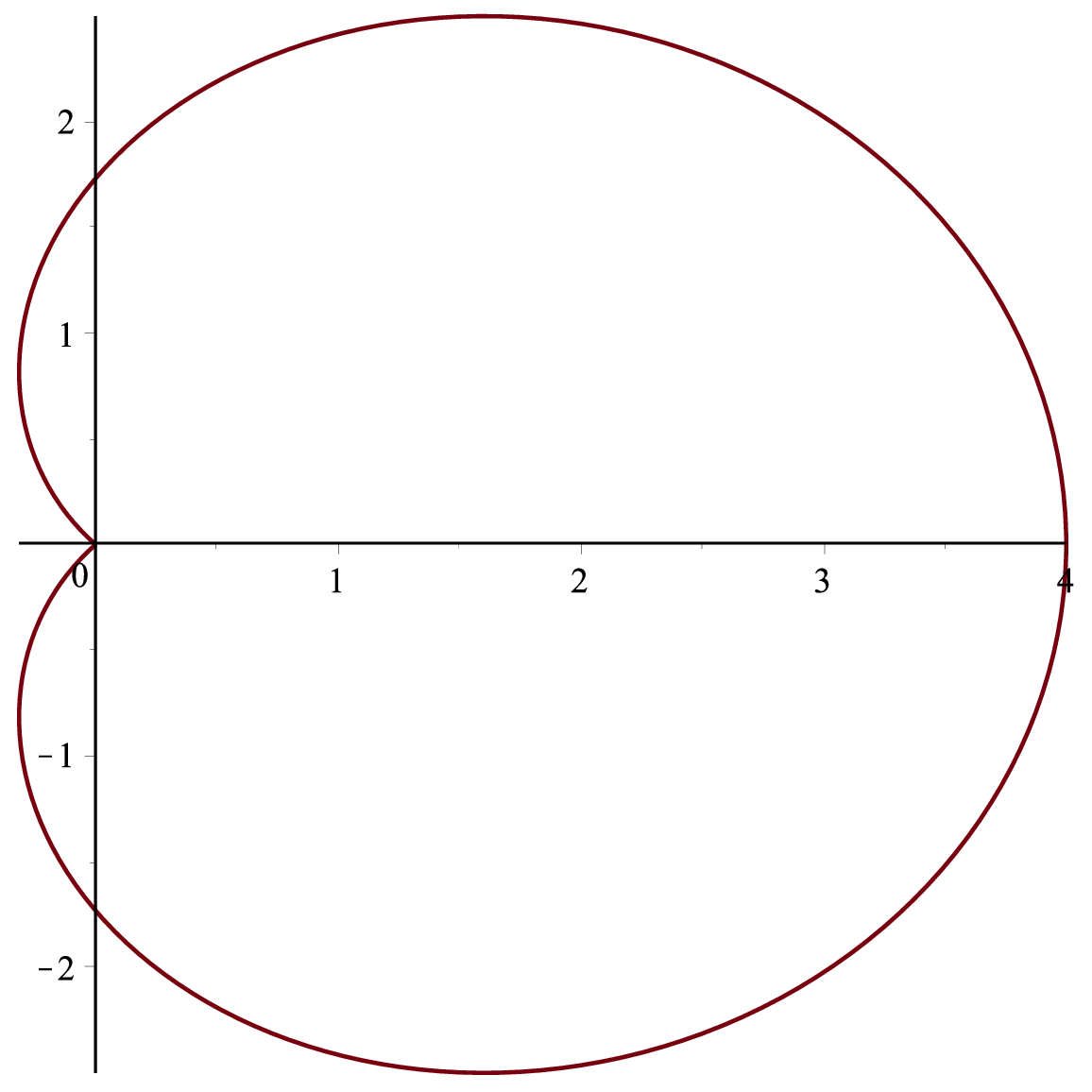} 
\end{minipage}
\hspace*{3cm}
\begin{minipage}{5cm}
   \includegraphics[width=5cm]{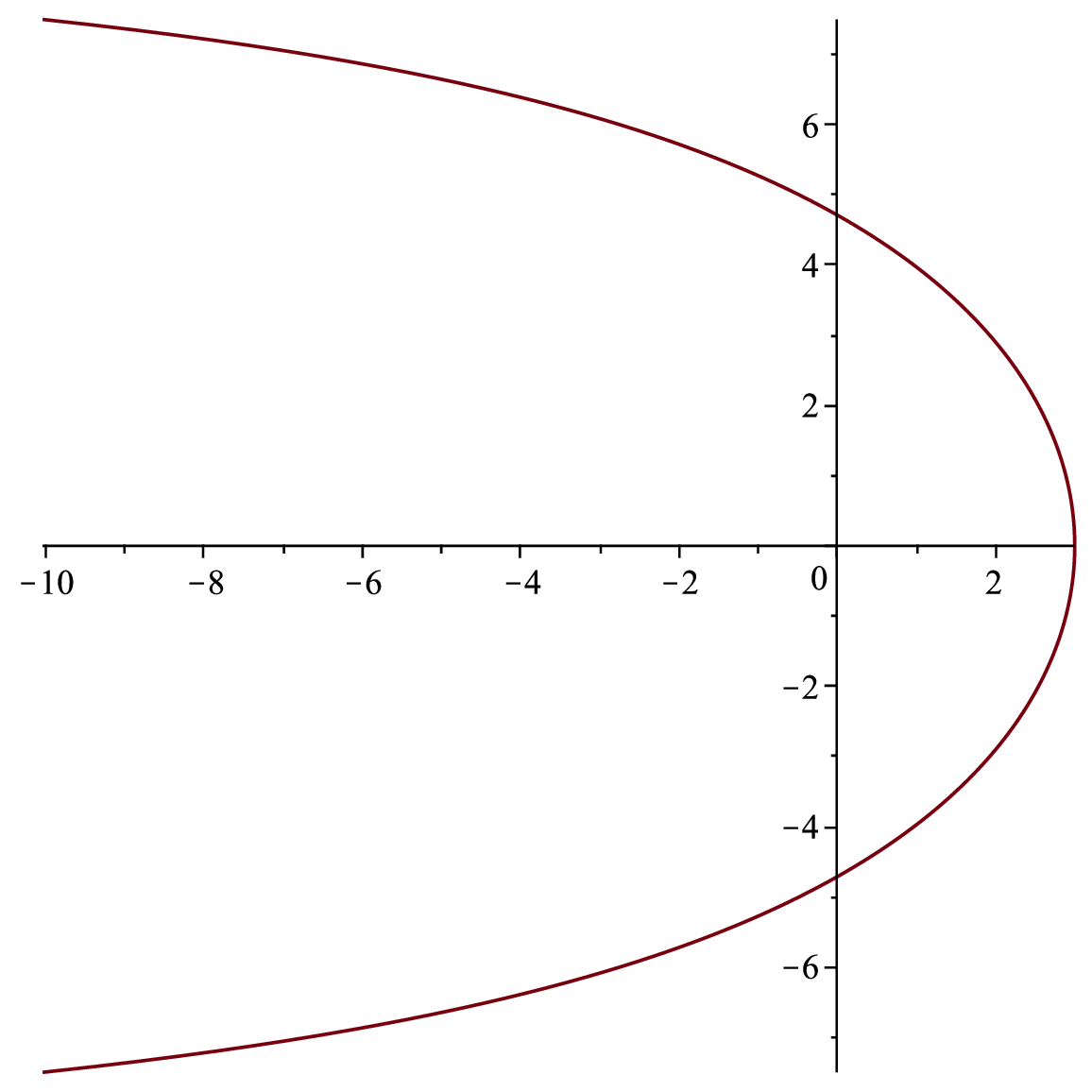} 
\end{minipage}
\end{center}
\caption{\small
Left: the steepest descent path of the integral in \eqref{eq:Mblea02} described by equation \eqref{eq:Mblea05}. Right: the steepest descent path of the integral in \eqref{eq:Mblea08}. In both cases we take  $\mu=3$. }
\label{fig:fig02}
\end{figure}

The saddle point $t_0$ follows from $\phi^\prime(t)=0$, where
\begin{equation}\label{eq:Mblea04}
\phi^\prime(t)=\frac{t-\mu-1}{t(t-1)}\quad \Longrightarrow\quad t_0=1+\mu.
\end{equation}
The path of steepest descent of the integral in \eqref{eq:Mblea02} through $t_0$ follows from the equation $\Im \phi(t)=0$. Using polar coordinates $t=r\cdot e^{i\theta}$ we find that it is given by
\begin{equation}\label{eq:Mblea05}
r=\frac{\sin((1+\mu)\theta/\mu)}{\sin(\theta/\mu)}, \quad -\frac{\mu}{1+\mu}\pi\le \theta\le \frac{\mu}{1+\mu}\pi.
 \end{equation}
In Figure~\ref{fig:fig02} (left) we show this path for $\mu=3$.

The standard saddle point method is not valid when $b\uparrow a$ and we use a uniform method transforming the integral in \eqref{eq:Mblea02} into the standard form \eqref{eq:appA02}. We use the transformation
\begin{equation}\label{eq:Mblea06}
\phi(t)-\phi(t_0)=\psi(s)-\psi(s_0),\quad \sign(t-t_0)=\sign(s-s_0),
\end{equation}
where
\begin{equation}\label{eq:Mblea07}
\psi(s)=s-\mu\ln s,\quad s_0=\mu;
\end{equation}
$s_0$ is the zero of $\psi^\prime(s)=(s-\mu)/s$.

We obtain
\begin{equation}\label{eq:Mblea08}
M(a,b,z)=\frac{\Gamma(b)}{\Gamma(a)}e^{bA} G_\lambda(b),\quad 
G_\lambda(b)=\frac{\Gamma(\lambda+1)}{2\pi i}
 \int _{-\infty}^{(0+)}e^{bs}s^{-\lambda-1}g(s)\,ds,
\end{equation}
where
\begin{equation}\label{eq:Mblea09}
  \frac{dt}{ds}=\frac{\psi^\prime(s)}{\phi^\prime(t)}=\frac{t(t-1)(s-s_0)}{s(t-t_0)},\quad 
  g(s)=\frac{s\,e^{zt}}{t(t-1)}\frac{dt}{ds}=e^{zt}\frac{s-s_0}{t-t_0},
\end{equation}
and
\begin{equation}\label{eq:Mblea10}
A=\phi(t_0)-\psi(s_0)=(1+\mu)\ln(1+\mu)-\mu.
\end{equation}
The saddle point contour of the integral in \eqref{eq:Mblea08} is the image of the contour in the $t$-plane described in \eqref{eq:Mblea05}. It runs through  $s_0=\mu$ and is defined by $\Im\psi(s)=0$. With polar coordinates $s=\rho\cdot e^{i\tau}$, we see that the contour is given by $\rho=\mu \tau/\sin\tau$, with $-\pi<\tau<\pi$. In Figure~\ref{fig:fig02} (right) we show this path for $\mu=3$.

Using the expansion given in \eqref{eq:appA18}, we obtain
\begin{equation}\label{eq:Mblea11}
M(a,b,z)\sim\frac{\Gamma(b)}{\Gamma(a)}e^{bA}  b^{\lambda}\sum_{n=0}^\infty (-1)^n
\frac{g_n(\mu)}{b^n},\quad a, b\to\infty, \quad a\ge b.
\end{equation}
To find $g_0(\mu)$ we evaluate (see the explanation as given for obtaining $f_0(\mu)$ in \eqref{eq:Mbgea12})
\begin{equation}\label{eq:Mblea12}
\left(\left.\frac{dt}{ds}\right\vert_{s=s_0}\right)^2=
\frac{\psi^{\prime\prime}(s_0)}{\phi^{\prime\prime}(t_0)}=1+\mu\quad
\Longrightarrow \quad g_0(\mu)=\frac{e^{z(1+\mu)}}{\sqrt{1+\mu}}.
\end{equation}

We take the coefficient $g_0(\mu)$ in front of the expansion and write 
\begin{equation}\label{eq:Mblea13}
M(a,b,z)\sim e^{bA}\frac{\Gamma(b)}{\Gamma(a)}b^{\lambda}g_0(\mu)\sum_{n=0}^\infty(-1)^n \frac{\wt{g}_n(\mu)}{b^n}, \quad
\wt{g}_n(\mu)=\frac{g_n(\mu)}{g_0(\mu)}, \quad a,b\to\infty,\quad a\ge b.
\end{equation}
We evaluate the  front factors by using the definition of $A$ in \eqref{eq:Mblea10} and  the scaled gamma functions defined in \eqref{eq:appB09}, and obtain
\begin{equation}\label{eq:Mblea14}
e^{bA}\frac{\Gamma(b)}{\Gamma(a)}b^{\lambda}g_0(\mu)=e^{z(1+\mu)}\frac{\Gamma^*(b)}{\Gamma^*(a)}.
\end{equation}
This gives the final result
\begin{equation}\label{eq:Mblea15}
M(a,b,z)\sim e^{z(1+\mu)}\frac{\Gamma^*(b)}{\Gamma^*(a)}\sum_{n=0}^\infty (-1)^n\frac{\wt{g}_n(\mu)}{b^n}, \quad a, b\to\infty,\quad b\le a.
\end{equation}
If we wish, we can expand the ratio of scaled gamma functions in front of this expansion  in powers of~$b^{-1}$ 
by  using $a=b(1+\mu)$.

The first few coefficients of this expansion are ${\wt g}_0(\mu)=1$,
\begin{equation}\label{eq:Mblea16}
\begin{array}{r@{\,}c@{\,}l}
{\wt g}_1(\mu)&=&\dsp{\frac{\mu\left(6z^2\nu^2+1\right)}{12(1+\mu)},\quad \nu=1+\mu }\\[8pt]
{\wt g}_2(\mu)&=&\dsp{\frac{\mu\left(36z^4\nu^4\mu+96z^3\nu^3(2\mu+1)+12z^2\nu^2(13\nu-1)-1\right)}{288(1+\mu)^2}.}
\end{array}
\end{equation}
The coefficients $g_n(\mu)$ are linear combinations of the derivatives of  $g(s)$ at the saddle point $s=\mu$ and  follow from the scheme given in Appendix~A.

The function  $g(s)$ defined in \eqref{eq:Mblea09} behaves as $g(s)=\bigO(se^{(1+\mu)\Re z})$ as $s,\mu\to\infty$, because $t_0=1+\mu$  and the path of integration in $t$-plane extends to $t$-values of $\bigO(\mu)$.  Although  the expansion in \eqref{eq:Mblea15} is scaled  by putting the exponential  $e^{(1+\mu)z}$ in front of the expansion, it is not possible to give a uniform bound for all $\mu\ge0$ of the iterates of $g(s)$ in the remainder. Therefore we have given the condition in \eqref{eq:Mblea01} on $\mu$ to be bounded. From the shown coefficients in \eqref{eq:Mblea16} we also see that $\mu$ should be bounded, except when $z=0$.

\section{\protectbold{U(a,b+1,z),\ b\ge a}}\label{sec:Ubgea}
For the $U$-function we consider $U(a,b+1,z)$ because this yields similar results as for $M(a,b,z)$. We have the special value $U(a,a+1,z)= z^{-a}$. 

In this section we use the notation and condition
\begin{equation}\label{eq:Ubgea01}
\lambda = b-a,\quad \mu=\frac{\lambda}{a}=\frac{b-a}{a},\quad \Re z>0,\quad \vert z\vert \le z_0,
\end{equation}
where $z_0$ is a fixed positive number.

We use the contour integral in \eqref{eq:appB06} and the Kummer relation for the $U$-function. This gives
\begin{equation}\label{eq:Ubgea02}
U(a,b+1,z)=\frac{z^{-b}\Gamma(\lambda+1)}{2\pi i}\int _{{-\infty}}^{{(0+)}} e^{zt}t^{-\lambda-1}(1-t)^{-a}\,dt.
\end{equation}
We write this in the form
\begin{equation}\label{eq:Ubgea03}
U(a,b+1,z)=\frac{z^{-b}\Gamma(\lambda+1)}{2\pi i}\int _{{-\infty}}^{{(0+)}}e^{a\phi(t)}e^{zt}\,\frac{dt}{t},
\end{equation}
where
\begin{equation}\label{eq:Ubgea04}
\phi(t)=-\ln(1-t)-\mu\ln t.
\end{equation}
The saddle point $t_0$ follows from 
\begin{equation}\label{eq:Ubgea05}
\phi^\prime(t)=\frac{(1+\mu)t-\mu}{t(1-t)}=0\quad \Longrightarrow \quad t_0=\frac{\mu}{1+\mu}.
\end{equation}

The saddle point contour is the curve through $t_0$ defined by $\Im\phi(t)=0$. We write $t=r\cdot e^{i\theta}$ and it follows that the contour is given by
\begin{equation}\label{eq:Ubgea06}
r=\frac{\sin(\mu\theta)}{\sin((1+\mu)\theta)},\quad -\frac{\pi}{1+\mu}<\theta<\frac{\pi}{1+\mu}.
\end{equation}

We use the transformation 
\begin{equation}\label{eq:Ubgea07}
\phi(t)-\phi(t_0)=\psi(s)-\psi(s_0), \quad \psi(s)=s-\mu\ln s,
\end{equation}
where $s_0=\mu$ is the zero of $\psi^\prime(s)$. 
This gives the representation
\begin{equation}\label{eq:Ubgea08}
U(a,b+1,z)=z^{-b}e^{aA} G_\lambda(a),\quad G_\lambda(a)=\frac{\Gamma(\lambda+1)}{2\pi i}\int _{{-\infty}}^{{(0+)}}e^{as}s^{-\lambda-1} p(s)\,ds,
\end{equation}
where
\begin{equation}\label{eq:Ubgea09}
p(s)=e^{zt}\frac{s}{t}\frac{dt}{ds}=e^{zt}\frac{(1-t)(s-\mu)}{(1+\mu)(t-t_0)}, \quad A=\phi(t_0)-\psi(s_0)=(1+\mu)\ln(1+\mu)-\mu.
\end{equation}
The saddle point contour in the $s$-plane is the same as the one for the integral in \eqref{eq:Mblea08}; see the right figure in Figure~\ref{fig:fig02}.

We have the expansion
\begin{equation}\label{eq:Ubgea10}
G_\lambda(a)\sim a^\lambda \sum_{n=0}^\infty (-1)^n\frac{p_n(\mu)}{a^n},\quad a,b\to\infty,\quad b\ge a.
\end{equation}
The first coefficient is
\begin{equation}\label{eq:Ubgea11}
p_0(\mu)=\frac{e^{z\mu/(1+\mu)}}{\sqrt{1+\mu}}.
\end{equation}

The first-order asymptotic approximation is
\begin{equation}\label{eq:Ubgea12}
U(a,b+1,z)\sim z^{-b}a^{b-a}e^{aA}p_0(\mu).
\end{equation}
Using the definition of $A(\mu)$ given in \eqref{eq:Ubgea09} 
 this becomes
\begin{equation}\label{eq:Ubgea13}
U(a,b+1,z)\sim z^{-b} a^{-a}b^b e^{a-b}  p_0(\mu).
\end{equation}
When $a=b$, that is, $\mu=0$, we obtain the value $z^{-a}$, which is the special value given in \eqref{eq:appB08}.

The full expansion can be written as
\begin{equation}\label{eq:Ubgea14}
U(a,b+1,z)\sim z^{-b} a^{-a}b^b e^{a-b} p_0(\mu)\sum_{n=0}^\infty (-1)^n\frac{\wt{p}_n(\mu)}{a^n},\quad a,b\to\infty,\quad b\ge a,
\end{equation}
where $\wt{p}_n(\mu)={p}_n(\mu)/p_0(\mu)$. We have $\wt{p}_0(\mu)=1$ and
\begin{equation}\label{eq:Ubgea15}
\begin{array}{r@{\,}c@{\,}l}
\wt{p}_1(\mu)&=&\dsp{\frac{\mu\left(\nu^2+6z(z-2-2\mu)\right)}{12\nu^3},\quad \nu=1+\mu,}\\[8pt]
\wt{p}_2(\mu)&=&\dsp{\frac{\mu\left(\mu\nu^4-24(\mu-12)\nu^3z +
 12(25\mu-36)\nu^2z^2-48(5\mu-2)\nu z^3+36\mu z^4 \right)}{288\nu^6}.}
\end{array}
\end{equation}

The function $p(s)$  defined in \eqref{eq:Ubgea09} behaves like $\bigO(s)$ as $s,\mu\to\infty$, the exponential function not being relevant in this case. By using the rational function representations as mentioned in Remark~\ref{rem:rem01}, we can find a uniform bound of the remainder in the expansion. The shown coefficients in \eqref{eq:Ubgea15} indicate that large values of $\mu$ are allowed.

\section{\protectbold{U(a,b+1,z),\ b\le a}}\label{sec:Ublea}
In this section we use the notation and conditions

\begin{equation}\label{eq:Ublea01}
\lambda = a-b,\quad \mu=\frac{\lambda}{b}=\frac{a-b}{b},\quad \mu\le \mu_0,\quad \Re z>0, \quad \vert z\vert \le z_0,
\end{equation}
where $\mu_0$ and $z_0$ are fixed positive numbers. The condition on $\mu$  means that $ b\ge a/(1+\mu_0)$, and that, say, $b=o(a)$ is not allowed.

We use the Kummer relation in \eqref{eq:appB07} and the integral representation in \eqref{eq:appB05}. This gives
\begin{equation}\label{eq:Ublea02}
U(a,b+1,z)=\frac{z^{-b}}{\Gamma(a-b)} \int_{0}^\infty e^{-zt}t^{a-b-1}(1+t)^{-a}\,dt,\quad \Re(a-b)>0,\quad \Re z>0,
\end{equation}
which we write in the form
\begin{equation}\label{eq:Ublea03}
U(a,b+1,z)=\frac{z^{-b}}{\Gamma(\lambda)} \int_{0}^\infty e^{-zt}e^{-b\phi(t)} \,\frac{dt}{t},
\end{equation}
where
\begin{equation}\label{eq:Ublea04}
\phi(t)=(1+\mu)\ln(1+t)-\mu\ln t.
\end{equation} 
We calculate the saddle point $t_0$: 
\begin{equation}\label{eq:Ublea05}
\phi^\prime(t)=\frac{t-\mu}{t(1+t)}=0\quad \Longrightarrow \quad t_0=\mu.
\end{equation}
We use the function $\psi(s)=s-\mu\ln s$ and transform
\begin{equation}\label{eq:Ublea06}
\phi(t)-\phi(t_0)=\psi(s)-\psi(s_0),\quad s_0=\mu,\quad \sign(t-t_0)=\sign(s-s_0),
\end{equation}
and write the result in the standard form
\begin{equation}\label{eq:Ublea07}
U(a,b+1,z)=z^{-b}e^{-bA}F_\lambda(b),\quad 
F_\lambda(b)=\frac{1}{\Gamma(\lambda)} \int_{0}^\infty e^{-bs}s^{\lambda-1}q(s)\,ds,
\end{equation}
where
\begin{equation}\label{eq:Ublea08}
q(s)= e^{-zt}\frac{s}{t} \frac{dt}{ds}=e^{-zt}\frac{(1+t)(s-\mu)}{t-\mu},\quad A=\phi(t_0)-\psi(s_0)=(1+\mu)\ln(1+\mu)-\mu.
\end{equation}

We have the expansion
\begin{equation}\label{eq:Ublea09}
F_\lambda(b)\sim b^{-\lambda}\sum_{n=0}^\infty \frac{q_n(\mu)}{b^n},\quad a,b\to\infty,\quad a\ge b.
\end{equation}

The first-order asymptotic approximation is
\begin{equation}\label{eq:Ublea10}
U(a,b+1,z)\sim z^{-b}b^{b-a}e^{-bA(\mu)}q_0(\mu),\quad q_0(\mu)=e^{-z \mu}\sqrt{1+\mu}.
\end{equation}
Using the definition of $A(\mu)$ given in \eqref{eq:Ublea09} 
 this becomes
\begin{equation}\label{eq:Ublea11}
U(a,b+1,z)\sim z^{-b} a^{-a}b^b e^{a-b}  q_0(\mu).
\end{equation}
When $a=b$, that is, $\mu=0$, we obtain the value $z^{-a}$, which is the special value given in \eqref{eq:appB08}.

The full expansion can be written as
\begin{equation}\label{eq:Ublea12}
U(a,b+1,z)\sim  z^{-b} a^{-a}b^b e^{a-b}  q_0(\mu)\sum_{n=0}^\infty \frac{\wt{q}_n(\mu)}{b^n},\quad a, b\to\infty,\quad b\le a,
\end{equation}
where $\wt{q}_n(\mu)={q}_n(\mu)/q_0(\mu)$. We have $\wt{q}_0(\mu)=1$ and
\begin{equation}\label{eq:Ublea13}
\begin{array}{r@{\,}c@{\,}l}
\wt{q}_1(\mu)&=&\dsp{\frac{\mu\left(6z^2\nu^2-12z\nu+1\right)}{12\nu},\quad \nu=\mu+1,}\\[8pt]
{\wt q}_2(\mu)&=&\dsp{\frac{\mu\left(36z^4\nu^4\mu-48z^3\nu^3(7\nu-5)+12z^2\nu^2(61\nu-25)-24z\nu(13\nu-1)+\mu\right)}{288\nu^2}.}
\end{array}
\end{equation}

Again, as in \S\ref{sec:Mblea}, we see that the coefficients grow with large values of $\mu$, and that we need to use the condition as shown in \eqref{eq:Ublea01}. Although the exponential function $e^{-zt}$ can be bounded uniformly  for $t\ge0$, this function has its influence in the $s$-variable. For large $t$ and  $s$, the transformation in \eqref{eq:Ublea06} takes the form $\ln t\sim s$, or $t\sim e^s$. Because for the evaluation of the coefficients we need values of the derivatives of the  function $q(s)$ at $s=\mu$, the exponential function has much influence on computing a uniform bound. When we take $z=0$ in the coefficients, we notice  the influence of the exponential function: the coefficients are bounded functions of $\mu$. Recall that $z=0$ is not allowed in this section.

\section{Numerical evaluations}\label{sec:num}
We give details on the numerical implementation of the expansions, and we  consider the case of \S\ref{sec:Mbgea} for $M(a,b,z)$, $b\ge a$.

The transformation in \eqref{eq:Mbgea05} can be written in the form
\begin{equation}\label{eq:num01}
\sum_{n=2}^\infty \frac{1}{n!}\phi^{(n)}(t_0)(t-t_0)^n=
\sum_{n=2}^\infty \frac{1}{n!}\psi^{(n)}(s_0)(s-s_0)^n,
\end{equation}
where the series converge in certain neighbourhoods of $t_0$ and $s_0$.  To invert the transformation near the saddle points, that is, to find $t$ when $s$ is given, we use the expansion $\dsp{t=t_0+\sum_{n=0}^\infty t_k(s-s_0)^n}$, and find $t_k$ by standard inversion methods for formal series. We have $t_1=\sqrt{\psi^{\prime\prime}(s_0)/\phi^{\prime\prime}(t_0)}=(1+\mu)^{-\frac32}$, where the square root has to be positive, in agreement with the condition $\sign(t-t_0)=\sign(s-s_0)$ imposed on the transformation in \eqref{eq:Mbgea05}.
The next terms are
\begin{equation}\label{eq:num02}
t_2=-\frac{\sqrt{\mu+1}\,(\mu-1)+\mu+1}{3\mu(\mu+1)^{\frac52}},\quad
t_3=\frac{\sqrt{\mu+1}\,(\mu^2-4\mu+8)+8\mu^2-8}{36\mu^2(\mu+1)^3}.
\end{equation}
These are analytic at  $\mu=0$, and we have
\begin{equation}\label{eq:num03}
t_2=-\tfrac12+\tfrac{25}{24}\mu+\bigO\left(\mu^2\right), \quad
t_3=\tfrac16-\tfrac{11}{24}\mu+\bigO\left(\mu^2\right),\quad \mu\to 0.
\end{equation}

The next step is to find the coefficients $a_n(\mu)$ in the expansion \eqref{eq:appA08}, with $f(s)$ defined in \eqref{eq:Mbgea08}, and finally we compute  the coefficients $f_n(\mu)$ by using the relations in \eqref{eq:appA11}. The first  scaled versions of these coefficients of the expansion in \eqref{eq:Mbgea15} are given in  \eqref{eq:Mbgea16}.

\renewcommand{\arraystretch}{1.2}
\begin{table}
\caption{
Relative errors in the computation of $M(a,b,z)$ for $b=1010.2$, $z=2.5$, several values of $a$  by using expansion\eqref{eq:Mbgea15} with terms up to $n=5$. The errors are computed by using the recurrence relation in \eqref{eq:num04}.\label{tab:tab01}}
$$
\begin{array}{rcccccc}
a\quad & n=0  & n=1 &  n=2 & n=3 &n=4 & n=5 \\
\hline
      5.1& \  0.27e-02& \  0.32e-04& \  0.83e-05& \  0.39e-07& \  0.13e-06& \  0.19e-08 \\
  205.1& \  0.15e-06& \  0.18e-08& \  0.12e-11& \  0.13e-13& \  0.00e-00& \  0.10e-14 \\
  405.1& \  0.96e-07& \  0.27e-08& \  0.25e-11& \  0.08e-14& \  0.00e-00& \  0.00e-00\\
  605.1& \  0.84e-06& \  0.23e-08& \  0.79e-11& \  0.17e-13& \  0.00e-00& \  0.20e-15 \\
  805.1& \  0.20e-05& \  0.94e-09& \  0.39e-11& \  0.28e-13& \  0.00e-00& \  0.00e-00\\
1005.1& \  0.31e-05& \  0.78e-08& \  0.26e-10& \  0.10e-12& \  0.25e-14& \  0.80e-15 \\
\hline
\end{array}
$$
\end{table}
\renewcommand{\arraystretch}{1.0}

For a numerical verification of the expansion we have used the expansion \eqref{eq:Mbgea15}  with terms up to $n=5$ and we have used a stable recursion relation (see \eqref{eq:appB09}) in the form
\begin{equation}\label{eq:num04}
\frac{zM(a+1,b+1,z)+bM(a,b,z)}{bM(a+1,b,z)}-1=0
\end{equation}
to verify the relative error in the approximations. In Table~\ref{tab:tab01} we show these errors for $b=1010.2$, $z=2.5$; $n=0,1,2,3,4,5$ means that we have used terms up to and including index $n$. We notice, for each $n$, a rather uniform error  for all values of $a$, except for $a=5.1$. Computations are done with Maple, with $Digits=16$.

\section{Concluding remarks}\label{sec:conclude}
In Section 2 and 4 we have given expansions in negative powers of $a$, although in both sections $b\ge a$. In Sections 3 and 5, the expansions are in negative powers of $b$, although $a\ge b$. For the asymptotics it is not relevant which parameter to choose, because both $a$ and $b$ are assumed to be large. In Sections 3 and 5 the representation of the coefficients is more attractive with negative powers of $b$ than with negative powers of $a$. This choice has no influence on whether or not we can take large values of $\mu$, which is only possible in Sections 2 and 4, where $b\ge a$. It appears that $b\ge a$ gives a better asymptotic condition for this type of asymptotic expansion for the Kummer functions. The starting point of these investigations was to obtain expansions valid for $a\sim b$, which always corresponds with $\mu\sim 0$, and it is an extra bonus when we have expansions that are valid for larger values of $\mu$ as well.

\section{Appendix~A: The vanishing saddle point}\label{sec:appA}
The asymptotic methods that we consider in this paper are for integrals of Laplace-type of the form
\begin{equation}\label{eq:appA01}
F_\lambda(w)=\frac1{\Gamma(\lambda)}\int_0^\infty
s^{\lambda-1}e^{-ws}f(s)\,ds,
\end{equation}
with $w$ as a large parameter. The method is also for loop integrals of the form
\begin{equation}\label{eq:appA02}
G_\lambda(w)=\frac{\Gamma(\lambda+1)}{2\pi i}\int_{-\infty}^{(0+)}
s^{-\lambda-1}e^{ws}f(s)\,ds,
\end{equation}
where the contour runs from $-\infty$ with $\phase\,s=-\pi$, encircles the origin in anti-clockwise direction, and returns to $-\infty$ with $\phase\,s=\pi$. The negative axis is a branch cut and we assume that $s^{-\lambda-1}$ has real values for $s>0$ (when $\lambda$ is real). In this paper we assume that $ w>0$ and $\lambda \ge0$.

When Watson's lemma is used for the integral in \eqref{eq:appA01}, with $w$ as the large parameter, the parameter $\lambda$ is assumed to be fixed. On the other hand, when, say $\lambda=\bigO(w)$, Watson's lemma cannot be used.
When $w$ and $\lambda$ are large, the  dominant part of the integral in \eqref{eq:appA01} is  
\begin{equation}\label{eq:appA03}
s^\lambda\,e^{-ws}=e^{-w\psi(s)},\quad \psi(s)=s-\mu\ln s,\quad \mu=\frac{\lambda}{w}.
\end{equation}
The function $\psi$ has  a 
saddle point at  $s=\mu$. When $w$ is large and $\lambda$ is fixed $\mu$ tends to zero, and the saddle point vanishes. 
When $\mu$ is bounded away from zero, we can transform the integral by using Laplace's method.
To describe an alternative method,  we summarise the treatment given in  \cite{Temme:1985:LTI}; see also   \cite[Chapter~25]{Temme:2015:AMI}, where the method is called {\em the vanishing saddle point}.

Consider \eqref{eq:appA01} and write $f(s)=\bigl(f(s)-f(\mu)\bigr)+f(\mu)$.
Then we have
\begin{equation}\label{eq:appA04}
\begin{array}{@{}r@{\;}c@{\;}l@{}}
F_\lambda(w)
&=&\dsp{w^{-\lambda}f(\mu)-
\frac1{w\Gamma(\lambda)}\int_0^\infty
\frac{f(s)-f(\mu)}{s-\mu}\,de^{-w\psi(s)}} \\[8pt]
&=&\dsp{w^{-\lambda}f(\mu)+\frac1{w\Gamma(\lambda)}\int_0^\infty 
s^{\lambda-1}e^{-ws}
f_1(s)\,ds,}
\end{array}
\end{equation}
where
\begin{equation}\label{eq:appA05}
f_1(s)=s\frac d{ds}\frac{f(s)-f(\mu)}{s-\mu}.
\end{equation}
Continuing this procedure we obtain for $K=0,1,2,\ldots$
\begin{equation}\label{eq:appA06}
\begin{array}{@{}r@{\;}c@{\;}l@{}}
\dsp{w^{\lambda}\,F_\lambda(w)}&=&\dsp{\sum_{k=0}^{K-
1}\frac{f_k(\mu)}{w^k}+
\frac1{w^K}E_K(w,\lambda),}\\[8pt]
\dsp{f_k(s)}&=&\dsp{s\frac d{ds}\frac{f_{k-1}(s)-f_{k-1}(\mu)}{s-\mu},\quad 
k=1,2,\ldots,\quad 
f_0(s)=f(s),} \\ [8pt]
\dsp{E_K(w,\lambda)}&=&\dsp{\frac1{\Gamma(\lambda)}\int_0^\infty s^{\lambda-1}
e^{-ws}f_K(s)\,ds.}
\end{array}
\end{equation}
Eventually we obtain the complete asymptotic expansion
\begin{equation}\label{eq:appA07}
F_\lambda(w)\sim w^{-\lambda}\sum_{n=0}^\infty
\frac{f_n(\mu)}{w^n},\quad w\to\infty.
\end{equation}

The coefficients $f_n(\mu)$ can be expressed in terms of the coefficients $a_n(\mu)$, which are defined by 
\begin{equation}\label{eq:appA08}
f(s)=\sum_{n=0}^\infty a_n(\mu)(s-\mu)^n.
\end{equation}
To verify this, we write
\begin{equation}\label{eq:appA09}
f_n(s)=\sum_{m=0}^\infty c_m^{(n)}(s-\mu)^m.
\end{equation}
Then $a_m(\mu)=c_m^{(0)}$, $f_n(\mu)=c_0^{(n)}$ and we have from \eqref{eq:appA06}
\begin{equation}\label{eq:appA10}
f_{n+1}(s)=\sum_{m=0}^\infty c_m^{(n+1)}(s-\mu)^m=s\sum_{m=1}^\infty c_m^{(n)}(m-1)(s-\mu)^{m-2}.
\end{equation}
This gives the recursion
\begin{equation}\label{eq:appA11}
c_{m}^{(n+1)}=mc_{m+1}^{(n)}+\mu(m+1)c_{m+2}^{(n)},\quad m,n=0,1,2,\ldots,
\end{equation}
and the few first relations are\footnote{The reviewer observed: {\em It seems that the numerical coefficients are the same as the sequence A269940 in the OEIS. It would be worth investigating this in the future.} See  also https://oeis.org/A269940 .}
\begin{equation}\label{eq:appA12}
\begin{array}{ll}
f_0(\mu)= a_0(\mu),\quad f_1(\mu)= \mu a_2(\mu),\quad f_2(\mu)= \mu\left(2a_3(\mu)+3\mu a_4(\mu)\right),\\[8pt]
f_3(\mu)= \mu\left(6a_4(\mu)+20\mu a_5(\mu)+15\mu^2a_6(\mu)\right),\\[8pt]
f_4(\mu)= \mu\left(24a_5(\mu)+130\mu a_6(\mu)+210\mu^2a_7(\mu)+105\mu^3a_8(\mu)\right).
\end{array}
\end{equation}
The functions $f_n(s)$ can be written as Cauchy-type integrals. Write $R_0(\sigma,s,\mu)=1/(\sigma-s)$. Then 
\begin{equation}\label{eq:appA13}
f_n(s)=\frac{1}{2\pi i}\int_{\calC} R_n(\sigma,s,\mu)f(\sigma)\,d\sigma,\quad R_{n+1}(\sigma,s,\mu)=\frac{-1}{\sigma-\mu}\frac{d}{d\sigma}\left(\sigma R_n(\sigma,s,\mu)\right),
\end{equation}
where $\calC$ is a simple closed contour in  the domain where $f(s)$ is analytic, and encircles the points $s$ and $\mu$. 
For large values of $s$ and $\mu$ it is not needed to take a large contour around the points $\sigma=s$ and $s=\mu$, because the contour $\calC$ can be split up into two circles around these points.

The next rational functions are
\begin{equation}\label{eq:appA14}
R_1(\sigma,s,\mu)=\frac{s}{(\sigma-\mu)(\sigma-s)^2},\quad 
R_2(\sigma,s,\mu)=\frac{s(\mu s + \mu \sigma - 2 \sigma^2)}{(\sigma - \mu)^3(\sigma-s)^3}.
\end{equation}

Under mild conditions on $a_n(\mu)$, that is,
on $f$, the expansion in \eqref{eq:appA07} is uniformly valid with respect to 
$\lambda\in
[0,\infty)$, and in a larger domain in the complex plane. The 
main
condition on $f$ is that its singularities are not too close to
the point $t=\mu$ and that $\vert f(s)\vert$ is bounded by an algebraic factor. 

Initially we have assumed for the integral in \eqref{eq:appA01}  that  $\lambda > 0$. However,  the reciprocal gamma function $1/\Gamma(\lambda)$ in front of the integral makes the integral regular when $\lambda\downarrow 0$. This can be seen by using integration by parts (writing $s^{\lambda-1}\,ds=(1/\lambda)\,d\left(s^\lambda\right)$), and in this way it can be shown that analytic continuation of $F_\lambda(w)$ of \eqref{eq:appA01} is possible into the domain $\Re\lambda\ge0$. We will see that the asymptotic expansion of $F_\lambda(w)$ allows taking $\lambda=0$. In fact the obtained expansion will be valid for $w\to\infty$, uniformly with respect to $\lambda \ge0$.

A similar integration by parts procedure gives the expansion of the loop integral in \eqref{eq:appA02}.
We use the integral of the reciprocal gamma function
\begin{equation}\label{eq:appA15}
\frac{w^\lambda}{\Gamma(\lambda+1)}=\frac{1}{2\pi i}\int_{-\infty}^{(0+)}
s^{-\lambda-1}e^{ws}\,ds,
\end{equation}
where the contour is a Hankel loop as in \eqref{eq:appA02}. Writing $\mu=\lambda/w$ and $g(s)=g(\mu)+\left(g(s)-g(\mu)\right)$, we obtain 
\begin{equation}\label{eq:appA16}
G_\lambda(w)= w^\lambda g(\mu)+\frac{1}{2\pi i w}\int_{-\infty}^{(0+)} \frac{g(s)-g(\mu)}{s-\mu}\frac{d}{ds}\left(e^{w\psi(s)}\right), \quad \psi(s)=s-\mu\ln s.
\end{equation}
Performing integration by parts, and repeating the procedure gives
\begin{equation}\label{eq:appA17}
\begin{array}{@{}r@{\;}c@{\;}l@{}}
\dsp{w^{-\lambda}\,G_\lambda(w)}&=&\dsp{\sum_{k=0}^{K-
1}(-1)^k\frac{g_k(\mu)}{w^k}+
\frac1{w^K}E_K(w,\lambda),}\\[8pt]
\dsp{E_K(w,\lambda)}&=&\dsp{\frac{\Gamma(\lambda+1)}{2\pi i}\int_{-\infty}^{(0+)}
s^{-\lambda-1}e^{ws}g_K(s)\,ds,}
\end{array}
\end{equation}
where the coefficients $g_k(\mu)$ can be obtained by the same recursive scheme as for $f_k(\mu)$ shown in \eqref{eq:appA06}. Eventually this gives the expansion
\begin{equation}\label{eq:appA18}
G_\lambda(w)\sim w^\lambda \sum_{k=0}^\infty (-1)^k\frac{g_k(\mu)}{w^k},\quad w\to \infty.
\end{equation}
Under conditions on $g(s)$, this expansion holds uniformly with respect to $\lambda\ge0$.

\section{Appendix~B}\label{sec:appB}

The defining power series is
\begin{equation}\label{eq:appB01}
M(a,b,z)=\sum_{n=0}^\infty \frac{(a)_n}{(b)_n}\frac{z^n}{n!}, \quad (a)_n=\frac{\Gamma(a+n)}{\Gamma(a)},
\end{equation}
with the usual condition that $b$ is not a nonpositive integer. The standard integral is
\begin{equation}\label{eq:appB02}
M(a,b,z)=\frac{\Gamma(b)}{\Gamma(a)\Gamma(b-a)}\int_0^1 e^{zt}t^{a-1}(1-t)^{b-a-1}\,dt,
\end{equation}
where $\Re a>0, \ \Re(b-a)>0$. A contour integral is
\begin{equation}\label{eq:appB03}
M(a,b,z)=\frac{\Gamma(b)\Gamma(1+a-b)}{\Gamma(a)}\frac{1}{2\pi i} \int _0^{(1+)}e^{zs} s^{a-1}(s-1)^{b-a-1}\,ds, \quad \Re a>0,
\end{equation}
where the contour starts at $s=0$, encircles the point $s=1$ in the anti-clockwise direction, and returns to $s=0$. Also,
\begin{equation}\label{eq:appB04}
M(a,b,z)=\frac{\Gamma(b)z^{1-b}}{2\pi i}\int_{\cal C} e^{zs}s^{-b}\left(1-1/s\right)^{-a}\,ds,
\end{equation}
where the contour ${\cal C}$ starts at $-\infty$, with $\phase\,s=-\pi$, encircles the points $0$ and $1$ in anti-clockwise direction, and returns to $-\infty$, where $\phase\,s=+\pi$. At the point where the contour crosses the interval $(1,\infty)$
the functions  $s^{-b}$ and $(1-1/s)^{-a}$ assume their principal values.

The standard integral  for  $U(a,b,z)$ is
\begin{equation}\label{eq:appB05}
U(a,b,z)=\frac{1}{\Gamma(a)}\int_0^\infty e^{-zt} t^{a-1}(1+t)^{b-a-1}\,dt,\quad \Re a>0, \quad  \Re z >0,
\end{equation}
and a loop integral is 
\begin{equation}\label{eq:appB06}
U(a,b,z)=\frac{\Gamma(1-a)}{2\pi i}\int _{{-\infty}}^{{(0+)}}e^{{zs}}s^{{a-1}}{(1-s)^{{b-a-1}}}ds,\quad \Re z>0,
\end{equation}
where $a\ne1,2,3,\ldots$. The contour cuts the real axis 
between $0$ and~$1$. At this point the fractional powers are determined by $\phase\,(1-s)=0$ and $\phase\,s=0$.

The Kummer relations are
\begin{equation}\label{eq:appB07}
M(a,b,z)=e^zM(b-a,b,-z),\quad U(a,b,z)=z^{1-b}U(a-b+1,2-b,z).
\end{equation}
Special values are
\begin{equation}\label{eq:appB08}
M(a,a,z)=e^z,\quad U(a,a+1,z)=z^{-a}.
\end{equation}
In numerical computations we have used the relation 
\begin{equation}\label{eq:appB09}
zM(a+1,b+1,z)+bM(a,b,z)=bM(a+1,b,z)
\end{equation}
to check the relative accuracy.

We use also the scaled gamma function
\begin{equation}\label{eq:appB10}
\Gamma^*(z)=e^z z^{-z}\sqrt{\frac{z}{2\pi}}\Gamma(z)\sim1+\frac{1}{12z}+\frac{1}{288 z^2}+\ldots,\quad z\to\infty.
\end{equation}

\section*{Disclosure statement}
No potential conflict of interest was reported by the author.

\section*{Funding}
This work was supported by the Spanish Ministerio de Ciencia, Innovaci\'on y Universidades under
Grants MTM2015-67142-P (MINECO/FEDER, UE) and\\
 PGC2018-098279-B-I00 (MCIU/AEI/FEDER, UE). \\

\section*{Acknowledgements}
The author is grateful to the reviewer for careful reading earlier versions of
the manuscript and for helpful comments that improved the article. \\
The author  thanks CWI, Amsterdam, for scientific support.

\end{document}